\documentclass[10pt]{amsart}
\usepackage{amsfonts,amssymb,amscd,amsmath,enumerate,verbatim,calc,amsthm}
\input xy
\xyoption{all}

\headsep=1cm \numberwithin{equation}{section}
\newtheorem{thm}{Theorem}[section]
\newtheorem{cor}[thm]{Corollary}
\newtheorem{lem}[thm]{Lemma}
\newtheorem{prop}[thm]{Proposition}
\newtheorem{defn}[thm]{Definition}
\newtheorem{exas}[thm]{Examples}
\newtheorem{exam}[thm]{Example}
\newtheorem{rem}[thm]{Remark}

\newcommand{\Hom}{\mbox{Hom}\,}
\newcommand{\Ext}{\mbox{Ext}\,}
\newcommand{\Tor}{\mbox{Tor}\,}
\newcommand{\Spec}{\mbox{Spec}\,}
\newcommand{\Max}{\mbox{Max}\,}
\newcommand{\Ker}{\mbox{Ker}\,}
\newcommand{\Ass}{\mbox{Ass}\,}

\newcommand{\Supp}{\mbox{Supp}\,}

\newcommand{\gr}{\mbox{grade}\,}
\newcommand{\Gid}{\mbox{Gid}\,}
\newcommand{\Gfd}{\mbox{Gfd}\,}
\newcommand{\Gpd}{\mbox{Gpd}\,}
\newcommand{\vdim}{\mbox{vdim}\,}
\newcommand{\depth}{\mbox{depth}\,}
\renewcommand{\dim}{\mbox{dim}\,}
\renewcommand{\Im}{\mbox{Im}\,}
\newcommand{\cd}{\mbox{cd}\,}

\newcommand{\h}{\mbox{ht}\,}
\newcommand{\im}{\mbox{Im}\,}
\newcommand{\E}{\mbox{E}}

\renewcommand{\H}{\mbox{H}}
\newcommand{\V}{\mbox{V}}

\newcommand{\fa}{\mathfrak{a}}
\newcommand{\fb}{\mathfrak{b}}

\newcommand{\fm}{\mathfrak{m}}
\newcommand{\fp}{\mathfrak{p}}
\newcommand{\fq}{\mathfrak{q}}

\newcommand{\ind}{\mbox {injdim}}
\begin{document}
\title[injective and gorenstein injective dimensions]
 { on injective and gorenstein injective dimensions of local cohomology modules}

\bibliographystyle{99}

     \author[M.R. Zargar]{Majid Rahro Zargar }
     \author[H. Zakeri]{Hossein Zakeri }

\address{M.R. Zargar and H. Zakeri, Faculty of mathematical sciences and computer, Kharazmi University, 599 Taleghani Avenue, Tehran 15618, Iran}
\email{zargar9077@gmail.com}
\email{zakeri@tmu.ac.ir}
\address {M.R. Zargar, School of Mathematics, Institute for Research in Fundamental Sciences (IPM), P.O. Box 19395--5746, Tehran, Iran.}




\begin{abstract} Let $(R,\fm)$ be a commutative Noetherian local ring and let $M$ be an $R$-module which is a relative Cohen-Macaulay with respect to a proper ideal $\fa$ of $R$ and set $n:=\h_{M}\fa$. We prove that $\ind M<\infty$ if and only if $\ind\H^{n}_\fa(M)<\infty$ and that $\ind\H^{n}_\fa(M)=\ind M-n$. We also prove that if $R$ has a dualizing complex and $\Gid_{R} M<\infty$, then $\Gid_{R}\H^{n}_\fa(M)<\infty$ and $\Gid_{R}\H^{n}_\fa(M)=\Gid_{R} M-n$. Moreover if $R$ and $M$ are Cohen-Macaulay, then it is proved that $\Gid_{R} M<\infty$ whenever $\Gid_{R}\H^{n}_\fa(M)<\infty$. Next, for a finitely generated $R$-module $M$ of dimension $d$, it is proved that if
$K_{\widehat M}$ is Cohen-Macaulay and $\Gid_{R}\H_{\fm}^{d}(M)<\infty$, then
$\Gid_{R}\H_{\fm}^{d}(M)=\depth R-d.$ The above results have consequences which improve some known results and provide characterizations of Gorenstein rings.

\end{abstract}

\maketitle

\section{introduction}

Throughout this paper, $R$ is a commutative Noetherian ring, $\fa$ is a proper ideal of $R$ and $M$ is an $R$-module. For a prime ideal $\fp$ of $R$, the residue class field $R_{\fp}/\fp R_{\fp}$ is denoted by $k(\fp)$. For each non-negative integer $i$, let $\H_{\fa}^{i}(M)$ denotes the $i$-th local cohomology module of $M$ with respect to $\fa$; see [1] for its definition and basic results. Also, we use $\ind_{R} M$ to denote the usual injective dimension of $M$. The notion of Gorenstein injective module was introduced by E.E. Enochs and O.M.G. Jenda in [4]. The class of Gorenstein injective modules is greater than the class of injective modules; but they are same classes whenever $R$ is a regular local ring. The Gorenstein injective dimension of $M$, which is denoted by $\Gid_{R}M$, is defined in terms of resolutions of Gorenstein injective modules. This notion has been used in [3, 15, 21] and has led to some interesting results. Notice that $\Gid_{R} M\leq\ind_{R} M$ and the equality holds if $\ind_{R} M<\infty$.

The principal aim of this paper is to study the injective (resp. Gorenstein injective) dimension of certain $R$-modules in terms of injective (resp. Gorenstein injective) dimension of its local cohomology modules at support $\fa$.

In this paper we will use the concept of relative Cohen-Macaulay modules which has been studied in [7] under the title of cohomologically complete intersections. The organization of this paper is as follows. In section 2, among other things, we prove, in 2.5, that if $M$ is relative Cohen-Macaulay with respect to $\fa$, then $\ind M$ and $\ind\H_{\fa}^{\h_{M}\fa}(M)$ are simultaneously finite and there is an equality $\ind\H_{\fa}^{\h_{M}\fa}(M)=\ind M-\h_{M}\fa.$ Then, as a corollary, we obtain a characterization of Gorenstein rings. Next, in 2.10, for all $n\geq0$ and any $\fp\in\Supp(M)$, we establish a comparison between the Bass numbers of $\H_{\fp R_{\fp}}^n(M_{\fp})$ and $\H_{\fm}^{n+\dim({R}/{\fp})}(M)$ whenever $(R,\fm)$ is a homomorphic image of a Gorenstein ring and $M$ is finitely generated.

In section 3, we first prove some basic properties about Gorenstein injective dimension of a module. In particular, Proposition 3.6 indicates that Gorenstein injective dimension is a refinement of the injective dimension. As a main result, in Theorem 3.8 we establish a Gorenstein injective version of 2.5. Indeed, it is proved that if, in addition to the hypothesis of 2.5, $R$ has a dualizing complex, then $\Gid_{R}M<\infty$ implies $\Gid_{R}\H_{\fa}^n(M)<\infty$ and the converse holds whenever $R$ and $M$ are Cohen-Macaulay. This theorem has consequences which recover some interesting results that have currently been appeared in the literature. As a first corollary of 3.8, we deduce that $\Gid_{R}\H_{\fm}^n(M)=\Gid_{R}M-n$, wherever $M$ is a Cohen-Macaulay module over the Cohen-Macaulay local ring $(R,\fm)$ and $\dim M=n$. This corollary improves the main result [15, Theorem 3.10](see the explanation which is offered before 3.9). As a second corollary, we obtain a characterization of Gorenstein local rings which recovers [21, Theorem 2.6]. As a main result, it has been shown in [15, Theorem 3.10] that if $R$ and $M$ are Cohen-Macaulay with $\dim M=d$ and $\Gid_{R}\H_{\fm}^d(M)<\infty$, then $\Gid_{R}\H_{\fm}^d(M)=\dim{R}-d$.
In 3.12, we will use the canonical module of a module to improve the above result without Cohen-Macaulay assumption on $R$ and $M$. This result provides some characterizations of Gorenstein local rings.
\section{local cohomology and injective dimension}

The starting point of this section is the next proposition, which plays essential role in the proof of Theorems 2.5 and 3.8.
\begin{prop} Let $n$ be a non-negative integer and let $N$ be an $\fa$-torsion $R$-module. Suppose that
$\emph\H^{i}_\fa(M)=0$ for all $i\neq n$. Then $$\emph\Ext^i_R(N,\emph\H^{n}_\fa(M))\cong \emph\Ext^{i+n}_R(N,M)$$ for all $i\geq 0.$

\end{prop}
\begin{proof} First we notice that $\Hom_{R}(N,M)=\Hom_{R}(N,\Gamma_{\fa}(M))$. Hence, in view of [14, Theorem 10.47], we have the Grothendieck third quadrant spectral sequence with $$E^{p,q}_{2}=\Ext_{R}^{p}(N, \H^{q}_{\fa}(M))\underset{p}\Longrightarrow \Ext_{R}^{p+q}(N,M).$$
Now, since $\H^{q}_\fa(M)=0$ for all $q\neq n$, $E^{p,q}_{2}=0$ for all $q\neq n$. Therefore, this spectral sequence collapses in
the column $q=n$; and hence one gets, for all $i\geq0$, the isomorphism
$$\Ext^i_R(N,\H^{n}_\fa(M))\cong\Ext^{i+n}_R(N,M),$$ as required.
\end{proof}

The following corollary, which is an immediate consequence of 2.1, determines the Bass numbers $\mu^{i}(\fp,\H^{n}_\fa(M)):=\vdim_{k(\fp)}\Ext^i_{R_{\fp}}(k(\fp),\H^{n}_{\fa R_{\fp}}(M_{\fp}))$ of the local cohomology module $\H^{n}_\fa(M)$.
\begin{cor} Let $n$ and $M$ be as in \emph{2.1.} Then, for all $\fp\in\V(\fa)$, $\mu^{i}(\fp,\emph{H}^{n}_\fa(M))=\mu^{i+n}(\fp,M)$ for each $i\geq0$.
\end{cor}

\begin{defn} \emph{We say that a finitely generated $R$-module $M$ is relative Cohen–
Macaulay with respect to $\fa$ if there is precisely one non-vanishing local cohomology module of $M$ with
respect to $\fa$. Clearly this is the case if and only if $\gr(\fa,M)=\cd(\fa,M)$, where $\cd(\fa,M)$ is the largest integer $i$ for which $\H_{\fa}^i(M)\neq0$ and $\gr(\fa,M)$ is the least integer $i$ such that $\H_{\fa}^i(M)\neq0$.}
\end{defn}
Observe that the above definition provides a generalization of the concept of Cohen-Macaulay modules. Also, notice that the notion of relative Cohen-Macaulay modules is connected with the notion of cohomologically complete intersection ideals which has been studied in [7] and has led to some interesting results. Furthermore, such modules have been studied in [8] over certain rings.

\begin{rem}\emph{Let $M$ be a relative Cohen-Macaulay module with respect to $\fa$ and let $\cd(\fa,M)=n$. Then, in view of  [1, Theorems 6.1.4, 4.2.1, 4.3.2], it is easy to see that $\Supp\H^{n}_{\fa}(M)=\Supp({M}/{\fa M})$ and $\h_{M}\fa=\gr(\fa,M)$, where
$\h_{M}\fa =\inf\{\ \dim_{R_{\fp}}M_{\fp} |~ \fp\in\Supp(M/\fa M) ~\}$.}
\end{rem}

The following theorem, which is one of the main results of this section, provides a comparison between the injective dimensions of a relative Cohen-Macaulay  module and its non-zero local cohomology module. Here we adopt the convention that the injective dimension of the zero module is to be taken as $-\infty$.

\begin{thm}Let $(R,\fm)$ be local and let $n$ be a non-negative integer such that $\emph{H}^{i}_\fa(M)=0$ for all $i\neq n$.
\begin{itemize}
\item[(i)]{ If \emph{$\ind M<\infty$}, then \emph{$\ind\H^{n}_\fa(M)<\infty$.}}
\item[(ii)]{The converse holds whenever $M$ is finitely generated}.
\end{itemize}
Furthermore, if $M$ is non-zero finitely generated, then \emph{$\ind\H^{n}_\fa(M)=\ind M-n$}.
\begin{proof}

(i). Let $s:=\ind M<\infty$. We may assume that $\H^{n}_\fa(M)\neq0$; and hence $s-n\geq0$. Therefore, in view of 2.2,  $\mu^{i+(s-n)}(\fp,\H^{n}_\fa(M))=0$ for all $\fp\in\Spec(R)$ and for all $i>0$; so that $\ind \H^{n}_\fa(M)\leq s-n$.

(ii). Suppose that $M$ is finitely generated. We first notice that $\H^{n}_{\fa}(M)=0$ if and only if $M=\fa M$; and this is the case if and only if $M=0$. Therefore we may assume that $\H^{n}_{\fa}(M)\neq0$. Suppose that $t:=\ind \H^{n}_\fa(M)<\infty$. Then there exists a prime ideal $\fq$ of $R$ such that $\mu^{t}(\fq,\H^{n}_\fa(M))\neq0$. Hence, by 2.2,
$\mu^{t+n}(\fq,M)\neq0$. Next we show that $\mu^{t+n+i}(\fp,M)=0$ for all $\fp\in\Spec(R)$ and for all $i>0$.
Assume the contrary. Then there exists a prime ideal $\fp$ of $R$ such that $\mu^{t+n+j}(\fp,M)\neq0$ for some $j>0$. Let $r:=\dim{R}/{\fp}$. Then, by [11, \S18, Lemma 4], we have $\mu^{t+n+j+r}(\fm,M)\neq0$. Hence, by 2.2, $\mu^{j+t+r}(\fm,\H^{n}_\fa(M))\neq0$ which is a contradiction in view of the choice of $t$.
Therefore,  $\ind M\leq t+n$. The final assertion is a consequence of (i) and (ii).
\end{proof}
\end{thm}

Next, we provide an example to show that if $R$ is non-local, then Theorem 2.5(ii) is no longer true. Also, in 3.11, we present two examples which show that 2.5(ii) and 2.5(i), respectively, are no longer true without the finiteness and the relative Cohen-Macaulayness assumptions on $M$.

\begin{exam}\emph{Suppose that $R$ is a non-local Artinian ring with $\ind R=\infty$. Let $\Max(R)=\{\fm_{1},...,\fm_{n}\}$. Then, in view of [17, Exercise 8.49], we have $R=\bigoplus_{\fm\in \Max(R)}\Gamma_{\fm}(R)$. Now, since the injective dimension of $R$ is infinite, there exists a maximal ideal $\fm_{t}$ of $R$ such that the injective dimension of $\Gamma_{\fm_{t}}(R)$ is infinite. Set $M:=\E_{R}({R}/{\fm_{s}})\bigoplus\Gamma_{\fm_{t}}(R)$, where $\fm_{s}\in\Max(R)$ with $\fm_{s}\neq\fm_{t}$. Then $M$ is a finitely generated $R$-module with infinite injective dimension and $\H^i_{\fm_{s}}(M)=0$ for all $i\neq 0$; but $\Gamma_{\fm_{s}}(M)$ is injective.}
\end{exam}

It is well-known that if $(R,\fm,k)$ is a $d$-dimensional local ring, then $R$ is Gorenstein if and only if $R$ is Cohen-Macaulay and $\H^{d}_{\fm}(R)\cong \E_R(k)$. The following corollary, which recovers this fact, is an immediate consequence of 2.5.
\begin{cor}Let $(R,\fm)$ be local and let $R$ be relative Cohen-Macaulay with respect to $\fa$. Then $R$ is Gorenstein if and only if \emph{$\ind\H_{\fa}^{\h_{R}{\fa}}(R)$} is finite.

In particular, if $\textbf{{x}}=x_{1},...,x_{n}$ is an $R$-sequence for some non-negative integer $n$, then  $R$ is Gorenstein if and only if \emph{$\ind\H^n_{(\textbf{{x}})}(R)$} is finite.
\end{cor}

The following proposition, which is needed in the proof of 3.8, provides an explicit minimal injective resolution for the non-zero local cohomology module of a relative Cohen-Macaulay module.

\begin{prop} Suppose that $M$ is relative Cohen-Macaulay with respect to $\fa$ and that \emph{$n =\cd(\fa,M)$}. Then\emph{$$0\longrightarrow \H^{n}_{\fa}(M)\longrightarrow\bigoplus_{\fp\in \V(\fa)}\mu^n(\fp,M)\E(R/\fp)\longrightarrow\bigoplus_{\fp\in \V(\fa)}\mu^{n+1}(\fp,M)\E(R/\fp)\longrightarrow\cdots$$} is a minimal injective resolution for\emph{ $ \H^{n}_{\fa}(M)$.} Furthermore, \emph{$\Ass_{R}\H^n_{\fa}(M)=\{\fp\in\V(\fa)|~\mu^n(\fp,M)\neq0\}.$}
\begin{proof}Let $$0\longrightarrow M\stackrel{d^{-1}} \longrightarrow \E^{0}(M)\stackrel{d^0}\longrightarrow \cdots\longrightarrow \E^{n-1}(M)\stackrel{d^{n-1}}\longrightarrow \E^{n}(M)\stackrel{d^n}\longrightarrow \E^{n+1}(M) \stackrel{d^{n+1}}\longrightarrow \cdots$$ be a minimal injective resolution for $M$. If there exists a prime ideal $\fp$ in $\V(\fa)$ with $\mu^{n-1}(\fp,M)\neq0$, then $\depth_{R_{\fp}}M_{\fp}\leq n-1$. On the other hand, since $\fp\in\Supp(M/\fa M)$, 2.4 implies that $\H^n_{\fa R_{\fp}}(M_{\fp})\neq0$. Therefore $n=\gr_{R_{\fp}}(\fa R_{\fp}, M_{\fp})\leq\depth_{R_{\fp}}M_{\fp}\leq n-1$ which is a contradiction.
It follows that $\Gamma_{\fa}(E^{n-1}(M))=0$; and hence we obtain the minimal injective resolution
$$0\longrightarrow \H^{n}_{\fa}(M)\longrightarrow\Gamma_{\fa}(\E^{n}(M))\longrightarrow\Gamma_{\fa}(\E^{n+1}(M))\longrightarrow\cdots$$ for $\H^{n}_{\fa}(M)$. Now, we may use this resolution to complete the proof.
\end{proof}
\end{prop}

The following elementary lemma, which is needed in the proof of the next theorem, can be proved by using a minimal free resolution for $M$ and the concept of localization.
\begin{lem}Let $(R,\fm,k)$ be local and let $M$ be finitely generated. Then, for any prime ideal $\fp$ of $R$, \emph{$\vdim_{k({\fp})}\Tor_i^{R_{\fp}}(k(\fp),M_{\fp})\leq\vdim_{k}\Tor_i^{R}(k, M)$ }for all $i\geq0$.
\end{lem}

The next theorem provides a comparison of Bass numbers of certain local cohomology modules.
\begin{thm}Suppose that $(R,\fm,k)$ is a local ring which is a homomorphic image of a Gorenstein local ring and that $M$ is finitely generated. Let $n,m$ be non-negative integers. Then \emph{$\mu^m(\fp,\H^{n}_{\fp}(M))\leq\mu^{m}(\fm,\H^{n+\dim R/\fp}_{\fm}(M))$} for all \emph{$\fp\in\Spec(R)$.}

\begin{proof} Let $(R',\fm')$ be a Gorenstein local ring of dimension $n'$ for which there exists a surjective ring homomorphism $f: R'\rightarrow R$.
Let $\fp$ be a prime ideal of $R$ and set $\fp'=f^{-1}(\fp)$. Now $R'_{\fp'}$ is a Gorenstein local ring and $\dim R'/\fp'=\dim R/\fp$. Since $R'$ is Gorenstein, we have $\dim R'_{\fp'}=\dim R'-\dim R'/\fp'$.  In view of [1, Exercise 11.3.1] there is, for each $j\in\mathbb{Z}$, an $R_{\fp}$-isomorphism $\Ext_{R'_{\fp'}}^j(M_{\fp},R'_{\fp'})\cong(\Ext_{R'}^j(M,R'))_{\fp}$. Also, by the Local Duality Theorem [1, Theorem 11.2.6], we have $\H^{n}_{\fp R_{\fp}}(M_{\fp})\cong\Hom_{R_{\fp}}\big(\Ext_{R'_{\fp'}}^{n'-n-t}(M_{\fp},R'_{\fp'}),\E_{R_{\fp}}(k({\fp}))\big)$ as $R_{\fp}$-modules, where $t:=\dim R/\fp$, and $\H^{n+t}_{\fm}(M)\cong \Hom_{R}(\Ext_{R'}^{n'-n-t}(M,R'),\E_{R}(k))$. It therefore follows that
\[\begin{array}{rl}
\Ext^m_{R_{\fp}}(k({\fp}),\H^{n}_{\fp R_{\fp}}(M_{\fp}))&\cong\Ext^m_{R_{\fp}}({k({\fp})},\Hom_{R_{\fp}}(\Ext_{R'_{\fp'}}^{n'-n-t}(M_{\fp},R'_{\fp'}),\E_{R_{\fp}}(k({\fp}))))\\
&\cong\Hom_{R_{\fp}}(\Tor_m^{R_{\fp}}(k({\fp}),\Ext_{R'_{\fp'}}^{n'-n-t}(M_{\fp},R'_{\fp'})),\E_{R_{\fp}}(k({\fp}))).
\end{array}\] and
\[\begin{array}{rl}
\Ext^m_{R}(k,\H^{n+t}_{\fm}(M))&\cong \Ext^m_{R}(k,\Hom_{R}(\Ext_{R'}^{n'-n-t}(M,R'),\E_{R}(k)))\\
&\cong\Hom_{R}(\Tor_m^{R}(k,\Ext_{R'}^{n'-n-t}(M,R')),\E_{R}(k)).
\end{array}\]
Since by 2.9

$\vdim_{k({\fp})}(\Tor_m^{R_{\fp}}(k({\fp}),\Ext_{R'_{\fp'}}^{n'-n-t}(M_{\fp},R'_{\fp'})))\leq \vdim_{k}(\Tor_m^{R}(k,\Ext_{R'}^{n'-n-t}(M,R')))$, one may use the above isomorphisms to complete the proof.
\end{proof}
\end{thm}
It is known as Bass's conjecture that if a local ring admits a finitely generated module of finite injective dimension, then it is a Cohen-Macaulay ring. For the proof of this fact the reader is referred to [12] and [13]. In the next corollary we shall use this fact and the concept of a generalized Cohen-Macaulay module. Recall that, over a local ring $(R,\fm)$, a finitely generated module of positive dimension is a generalized Cohen-Macaulay module if $\H_{\fm}^{i}(M)$ is finitely generated for all $0\leq i<\dim M$.
\begin{cor}Let the situation be as in \emph{2.10}. Then the following statements hold.
 \begin{itemize}
\item[(i)]{\emph{$\ind_{R_{\fp}}\H^{n}_{\fp R_{\fp}}(M_{\fp})\leq\ind_{R}\H^{n+\dim R/\fp}_{\fm}(M)$} for all prime ideals $\fp$ of $R$ and for any $n\geq0$.}
\item[(ii)]{If $M$ is generalized Cohen-Macaulay with dimension $d$ such that \emph{$\H^{d}_{\fm}(M)$} is injective, then $M_{\fp}$ is Gorenstein, in the sense of
    \emph{[18]}, for all \emph{$\fp\in\Supp(M)\setminus\{\fm\}$.}}
\end{itemize}
\begin{proof}(i) is clear by 2.10.

(ii) Let $\fp\in\Supp(M)\setminus\{\fm\}$. By [1, Exercise 9.5.7], $M_{\fp}$ is Cohen-Macaulay and $\dim M_{\fp}+\dim R/\fp =\dim M$. Hence, in view of (i) and 2.5, we have $\ind M_{\fp}=\dim M_{\fp}$. Therefore, by [18, Theorem 3.11], [2, Theorem 3.1.17] and Bass's conjecture,  $M_{\fp}$ is Gorenstein.
\end{proof}
\end{cor}


\section{local cohomology and gorenstein injective dimension}
We first recall some definitions that we will use in this section.
\begin{defn}\emph{Following [4], an $R$-module $M$ is said to be Gorenstein injective if there exists
a $\Hom({Inj},-)$ exact exact sequence
$$\cdots\rightarrow E_1\rightarrow E_0 \rightarrowE^0\rightarrow E^1\rightarrow\cdots $$
of injective $R$-modules such that $M= \Ker(E^0\rightarrow E^1)$. We say that an exact sequence
$$0\rightarrow M\rightarrow G^0 \rightarrow G^1\rightarrow G^2\rightarrow\cdots$$
of $R$-modules and $R$-homomorphisms is a Gorenstein injective resolution for $M$, if each $G^i$ is Gorenstein injective. We
say that $\Gid_R M\leq n$ if and only if $M$ has a Gorenstein injective resolution of length $n$. If there
is no shorter resolution, we set $\Gid_R M = n$. Dually, an $R$-module $M$ is said to be Gorenstein projective if there
is a $\Hom(-,{Proj})$ exact exact sequence
$$\cdots\rightarrow P_1\rightarrow  P_0 \rightarrow P^0 \rightarrow  P^1 \rightarrow \cdots$$
of projective $R$-modules such that $M = \Ker(P^0 \rightarrow P^1)$. Similarly, one can define the Gorenstein projective
dimension, $\Gpd_R M$, of $M$.}
\end{defn}

\begin{defn}\emph{For a local ring $R$ admitting the dualizing complex $D_{R}$, we denote by $K_{M}$ the canonical module of an $R$-module $M$, which is defined to be $$K_{M}=\H^{d-n}(\textbf{R}\Hom_{R}(M,D_{R})),$$
where $d=\dim R$ and $n=\dim M$. Note that if $R$ is Cohen-Macaulay, then $K_{R}$ coincides with the classical definition of the canonical module of $R$ which is denoted by $\omega_{R}$.}
\end{defn}

\begin{defn}\emph{ Let $R$ be a Cohen-Macaulay local ring of Krull dimension $d$ which admits a
canonical module $\omega_{R}$. Following [4], let $\mathcal{J}_{0}(R)$ be the class
of $R$-modules $M$ which satisfies the following conditions.}
\begin{itemize}
\item[(i)]{ \emph{$\Ext_{R}^i(\omega_{R},M)=0$} , for all $i>0$.}
\item[(ii)]{ \emph{$\Tor^{R}_{i}(\omega_{R},\Hom_{R}(\omega_{R},M))=0$}, for all $i>0$.}
\item[(iii)]{ The natural map \emph{$\omega_{R}\otimes_{R}\Hom_{R}(\omega_{R},M)\rightarrow M$} is an isomorphism.}
\end{itemize}
This class of $R$-modules is called the Bass class.
\end{defn}

\begin{defn}\emph{Following [19], let $\fa$ and $\fb$ be ideals of $R$. We set
$$W(\fa, \fb) = \{\ \fp\in\Spec(R) ~| ~\fa^n\subseteq \fp + \fb  ~ \text{for some integer}  ~n > 0\}.$$
For an $R$-module $M$, $\Gamma_{\fa,\fb}(M)$ denotes a submodule of $M$ consisting of all elements
of $M$ with support in $W(\fa, \fb)$, that is,
$\Gamma_{\fa,\fb}(M) = \{\ x \in M ~|~\Supp(Rx)\subseteq W(\fa, \fb)\}$.}
\end{defn}

The following lemma has been proved in [10, Lemma 4.2] and is of assistance in the proof of Proposition 3.6.
\begin{lem}Let $(R,\fm,k)$ be local. Then \emph{$$\Ext^{i}_{R}(\E_{R}(k),M)\cong \Ext^{i}_{{R}}(\E(k), M\otimes_{R}\hat R)\cong\Ext^{i}_{\hat{R}}(\E_{\hat R}(k), M\otimes_{R}\hat R)$$} for all $i\geq0$.
\end{lem}

Let $(R,\fm,k)$ be local and let $M$ be a non-zero non-injective $R$-module  of finite Gorenstein injective dimension. It was shown in [5, Corollary 4.4] that if $\Ext^{i}_{R}(\E,M)=0$ for all $i\geq0$ and all indecomposable injective $R$-modules $\E\neq \E_{R}(k)$, then $\Gid_{R}M=\sup\{ ~i ~|~ \Ext^{i}_{R}(\E_{R}(k),M)\neq0 \}.$ Our next proposition, which is concerned with this result, indicates that Gorenstein injective dimension is a refinement of the injective dimension. However we will use 3.6 and 3.7 to prove the main theorem 3.8.

\begin{prop}Let $(R,\fm,k)$ be local and let $M$ be non-zero with $\emph{\Gid}_{R}M<\infty$. If either $M$ is finitely generated or Artinian, then \emph{$$\Gid_{R}M=\sup\{ ~i ~|~ \Ext^{i}_{R}(\E_{R}(k),M)\neq0 \}.$$}
\end{prop}
\begin{proof} First assume that $M$ is finitely generated. Then, by [3, Theorem 3.24], $\Gid_{R}M=\Gid_{\widehat R}\widehat M$. Now, since $\widehat R$ is complete and $\widehat M$ is finitely generated as an $\widehat R$-module, the proof of [6, Proposition 2.2] in conjunction with [5, Corollary 4.4 ] implies that $$\Gid_{\widehat R}\widehat M=\sup\{~ i~ |~ \Ext^{i}_{\widehat R}(\E_{\widehat R}(\widehat R/\widehat\fm),\widehat M)\neq0 \}.$$ Therefore we can use 3.5 to complete the proof.

 Next, we consider the case where  $M$ is Artinian. By [15, Lemma 3.6 ] and [1, Exercise 8.2.4], $\Gid_{R}(M)=\Gid_{\widehat R}(M)$ and, by [3, Theorem 4.25], $\Gfd_{\widehat R}(\Hom_{\widehat R}(M, \E_{\widehat R}(\widehat R/\widehat\fm)))=\Gid_{\widehat R}(M)$, where, for an $R$-module $X$, $\Gfd_{R}(X)$, denotes the Gorenstein flat dimension of $X$. Therefore, since $\Hom_{\widehat R}(M, \E_{\widehat R}(\widehat R/\widehat\fm))$ is finitely generated as an $\widehat R$-module, in view of [3, Theorem 4.24] and [3. Theorem 1.10] we have the first equality in the next display
\[\begin{array}{rl}
\Gfd_{\widehat R}(\Hom_{\widehat R}(M, \E_{\widehat R}(\widehat R/\widehat\fm)))&=\sup\{\ ~ i |~ \Ext_{\widehat R}^{i}(\Hom_{\widehat R}(M, \E_{\widehat R}(\widehat R/\widehat\fm)),\widehat R)\neq0 \}\\
&=\sup\{\ ~ i |~ \Ext_{R}^{i}(\E_{R}(k),M)\neq0 \}.
\end{array}\]
The last equality follows from [10, Theorem 4.3], because $\E_{R}(k)$ and $M$ are Artinian.
\end{proof}

\begin{lem}Let $(R,\fm)$ be a Cohen-Macaulay local ring and let $M$ be finitely generated. Suppose that $x\in\fm$ is both $R$-regular and $M$-regular. Then the following statements are equivalent.
\begin{itemize}
\item[(i)]{$\emph{\Gid}_{R}M<\infty$.}
\item[(ii)]{$\emph{\Gid}_{R/x R}M/xM<\infty$}.
\end{itemize}
Furthermore, \emph{$\Gid_{R/x R}M/xM=\Gid_{R}M -1$}.
\end{lem}
\begin{proof}In view of [3, Theorem 3.24], we can assume that $R$ is complete; and hence it admits a canonical module $\omega_{R}$.

(i)$\Rightarrow$(ii). It follows from [3, Proposition 3.9] that $\Gid_{R}M/xM<\infty$. Therefore we can use [4, Proposition 10.4.22], [11, p.140,lemma 2] and [2, Theorem 3.3.5], to see that $\Gid_{R/x R}M/xM<\infty$.

(ii)$\Rightarrow$(i). By [4, Proposition 10.4.23], $M/xM\in\mathcal{J}_{0}(R/xR)$. Since $\omega_{R/xR}\cong\omega_{R}/x\omega_{R}$, in view of [11, p.140, lemma 2] we have $\Ext_{R}^i(\omega_{R},M/xM)=\Tor_i^{R}(\omega_{R},\Hom_{R}(\omega_{R}, M/xM))=0$ for all $i>0$ and $\omega_{R}/x\omega_{R}\otimes_{R/xR}\Hom_{R}(\omega_{R},M/xM)\cong M/xM$. Now, using the exact sequence $\Ext_{R}^i(\omega_{R}, M)\stackrel{x}\longrightarrow \Ext_{R}^i(\omega_{R}, M)\longrightarrow\Ext_{R}^i(\omega_{R}, M/xM)$ and Nakayama's lemma, we deduce that $\Ext_{R}^i(\omega_{R}, M)=0$ for all $i>0$. Thus we have the exact sequence
\begin{equation}
 0\longrightarrow\Hom_{R}(\omega_{R}, M)\stackrel{x}\longrightarrow \Hom_{R}(\omega_{R}, M) \longrightarrow\Hom_{R}(\omega_{R}, M/xM) \longrightarrow 0.
\end{equation}
 Now, we may use (3.1) and similar arguments as above to see that $\Tor_{R}^i(\omega_{R},\Hom_{R}(\omega_{R},M))=0$ for all $i>0$. Also, in view of [2, Lemma 3.3.2], we can see that $\omega_{R}\otimes_{R}\Hom_{R}(\omega_{R},M)\cong M$. Therefore by [4, Proposition 10.4.23], $\Gid_{R}M<\infty$. The final assertion is an immediate consequence of [3, Theorem 3.24].
\end{proof}
Theorem 3.8, which is a Gorenstein injective version of 2.5, is one of the main results of this section. As we will see, this theorem has consequences which recover some interesting results that have currently been appeared in the literature.
Here we adopt the convention that the Gorenstein injective dimension of the zero module is to be taken as $-\infty$.

\begin{thm}Suppose that the local ring $(R,\fm)$ has a dualizing complex and let $n$ be a non-negative integer such that $\emph{H}^{i}_\fa(M)=0$ for all $i\neq n$.
\begin{itemize}
\item[(i)]{ If $\emph{\Gid}_{R}M<\infty$, then $\emph{\Gid}_{R}\emph{H}^{n}_\fa(M)<\infty$}.
\item[(ii)]{The converse holds whenever $R$ and $M$ are Cohen-Macaulay}.
\end{itemize}
Furthermore, if $M$ is non-zero finitely generated with finite Gorenstein injective dimension, then $\emph{\Gid}_{R}\emph{H}^{n}_\fa(M)=\emph{\Gid}_{R}M-n.$

\end{thm}
\begin{proof}(i) Notice that if $\H^n_{\fa}(M)=0$, then there is nothing to prove. So, we may assume that $\H^n_{\fa}(M)\neq0$. Hence, by [20, Lemma 1.1], we have $n\leq d$, where $d=\Gid_{R} M$. Let
$$0\longrightarrow M\stackrel{d^{-1}} \longrightarrow G^{0}\stackrel{d^0}\longrightarrow G^{1}\stackrel{d^1}\longrightarrow \cdots\longrightarrow G^{n-1}\stackrel{d^{n-1}}\longrightarrow G^{n}\stackrel{d^n}\longrightarrow G^{n+1} \stackrel{d^{n+1}}\longrightarrow \cdots\longrightarrow G^{d-1}\stackrel{d^{d-1}}\longrightarrow G^{d}\longrightarrow  0$$
be a Gorenstein injective resolution for $M$. By applying the functor $\Gamma_{\fa}(-)$ on this exact sequence, we obtain the complex $$0\longrightarrow \Gamma_{\fa}(M) \stackrel{\Gamma_{\fa}(d^{-1})}\longrightarrow \Gamma_{\fa}(G^{0})\stackrel{\Gamma_{\fa}(d^{0})}\longrightarrow\Gamma_{\fa}(G^{1})\stackrel{\Gamma_{\fa}(d^{1})}\longrightarrow \cdots\longrightarrow \Gamma_{\fa}(G^{n-1})$$  $$\stackrel{\Gamma_{\fa}(d^{n-1})}\longrightarrow \Gamma_{\fa}(G^{n})\stackrel{\Gamma_{\fa}(d^{n})}
\longrightarrow \Gamma_{\fa}(G^{n+1}) \stackrel{\Gamma_{\fa}(d^{n+1})}\longrightarrow \cdots
\longrightarrow \Gamma_{\fa}(G^{d-1})\stackrel{\Gamma_{\fa}(d^{d-1})}\longrightarrow \Gamma_{\fa}(G^{d})\longrightarrow 0$$ in which, by [15, Theorem 3.2], $\Gamma_{\fa}(G^i)$ is Gorenstein injective for all $0\leq i\leq d$. If $n=0$ the result is clear. So suppose that $n>0$.  Now, since by [20, Lemma 1.1] each $G^i$ is $\Gamma_{\fa}$-acyclic for all $i$, we may use [1, Exercise 4.1.2 ] in conjunction with our assumption on local cohomology modules of $M$ to obtain the following two exact sequences$$0\longrightarrow \Gamma_{\fa}(M) \stackrel{\Gamma_{\fa}(d^{-1})}\longrightarrow \Gamma_{\fa}(G^{0})\stackrel{\Gamma_{\fa}(d^{0})}\longrightarrow \cdots\longrightarrow \Gamma_{\fa}(G^{n-1})\stackrel{\Gamma_{\fa}(d^{n-1})}\longrightarrow\Gamma_{\fa}(G^{n})\longrightarrow \frac{\Gamma_{\fa}(G^{n})}{\Im\Gamma_{\fa}(d^{n-1})}\longrightarrow 0$$ and

$$0\longrightarrow\Im\Gamma_{\fa}(d^{n})=\Ker\Gamma_{\fa}(d^{n+1})\hookrightarrow \Gamma_{\fa}(G^{n+1}) \longrightarrow \cdots\longrightarrow\Gamma_{\fa}(G^{d-1})\longrightarrow \Gamma_{\fa}(G^{d})\longrightarrow 0.$$
 But, by assumption, $\Gamma_{\fa}(M)=0$. Therefore, by using  the first above exact sequence and [4, Theorem 10.1.4], we see that $\frac{\Gamma_{\fa}(G^{n})}{\Im\Gamma_{\fa}(d^{n-1})}$ is Gorenstein injective. Notice that  $\H^{n}_{\fa}(M)=\frac{\ker{\Gamma{\fa}(d^{n})}}{\im\Gamma_{\fa}(d^{n-1})}$. Therefore, patching the second above long exact sequence together with the exact sequence $$0\longrightarrow \H^{n}_{\fa}(M)\longrightarrow\frac{\Gamma_{\fa}(G^{n})}{\Im\Gamma_{\fa}(d^{n-1})}\longrightarrow \frac{\Gamma_{\fa}(G^{n})}{\ker\Gamma_{\fa}(d^{n})}\longrightarrow 0,$$ yields the following long exact sequence $$0\longrightarrow \H^{n}_{\fa}(M)\longrightarrow\frac{\Gamma_{\fa}(G^{n})}{\Im\Gamma_{\fa}(d^{n-1})}\longrightarrow \Gamma_{\fa}(G^{n+1}) \longrightarrow \cdots
\longrightarrow \Gamma_{\fa}(G^{d-1})\longrightarrow \Gamma_{\fa}(G^{d})\longrightarrow 0.$$ Hence, $\Gid_{R}\H^{n}_\fa(M)\leq\Gid_{R}M-n$.

(ii) Suppose that $R$ and $M$ are Cohen-Macaulay. Since $\H_{\fm}^0(\E(R/\fm))=\E(R/\fm)$ and, for any non-maximal prime ideal $\fp$ of $R$, $\H_{\fm}^0(\E(R/\fp))=0$, we may apply 2.8 to see that $\H_{\fm}^i(\H_{\fa}^n(M))=\H_{\fm}^{n+i}(M)$ for all $i\geq0$. Therefore, we can use the Cohen-Macaulayness of $M$ to deduce that
\[ \H_{\fm}^i(\H_{\fa}^n(M))=\begin{cases}
       0 & \text{if $i\neq\dim M/\fa M$}\\
       \H^d_{\fm}(M) & \text{if $i=\dim M/\fa M,$}
       \end{cases} \]
where $d=\dim M$. Hence, by using part(i) for $\H_{\fa}^n(M)$, we have $\Gid_{R}\H^{d}_\fm(M)<\infty$. Now, we proceed by induction on $d$ to show that $\Gid_{R}M$ is finite. The case $d=0$ is clear. Let $d>0$ and assume that the result has been proved for $d-1$. Suppose that $x\in\fm$ is both $R$-regular and $M$-regular. Then one can use the induced exact sequence $$0\longrightarrow {\H}^{d-1}_\fm(M/xM)\longrightarrow {\H}^{d}_\fm(M) \longrightarrow {\H}^{d}_\fm(M) \longrightarrow 0$$  and [3, Proposition 3.9]  to see that $\Gid_{R}({\H}^{d-1}_\fm(M/xM))$ is finite. Hence, by inductive hypothesis, $\Gid_{R}M/xM$ is finite. Therefore, since, in view of [9, Corollary 6.2], $R$ admits a canonical module, one can use the same argument as in the proof of 3.7(i)$\Rightarrow$(ii) to deduce that $\Gid_{R/xR}M/xM<\infty$. Therefore $\Gid_{R}M$ is finite by 3.7. Now the result follows by induction.

For the final assertion, let $M$ be non-zero finitely generated with $\Gid_{R} M=s<\infty$. Then, by part(i), we have $\Gid_{R}\H^{n}_\fa(M)\leq s-n$. If $\Gid_{R}\H^{n}_\fa(M)< s-n$, Then, in view of [3, Theorem 3.6], we deduce that $\Ext_{R}^{s-n}(\E(k),\H^{n}_{\fa}{(M)})=0$. Hence, by Proposition 2.1, $\Ext_{R}^{s}(\E(k),M)=0$ which is a contradiction by 3.6. Therefore, $\Gid_{R}\H^{n}_\fa(M)=\Gid_{R}M-n.$
\end{proof}

Let $(R,\fm)$ be a local ring. As a main theorem, it was proved in [15, Theorem 3.10] that if $R$ and $M$ are Cohen-Macaulay with $\dim M=n$ and $\Gid_{R}{\H}^{n}_\fm(M)<\infty$, then $\Gid_{R}{\H}^{n}_\fm(M)=\dim R-n$. Notice that if $\Gid_{R}{\H}^{n}_\fm(M)<\infty$, then, in view of [15, Lemma 3.6], 3.8(ii) and [3, Theorem 3.24], we have $\depth R=\Gid_{R}M$. Therefore, the next corollary, which is established without the assumption that $\Gid_{R}{\H}^{n}_\fm(M)<\infty$, recovers [15, Theorem 3.10]. Another improvement of the above result will be given in 3.12.
\begin{cor}Let $(R,\fm)$ be a Cohen-Macaulay local ring and let $M$ be Cohen-Macaulay of dimension $n$. Then $\emph{\Gid}_{R}\emph{H}^{n}_\fm(M)=\emph{\Gid}_{R}M-n$.

\begin{proof} First we notice that $M \otimes_R\widehat R$ is a Cohen-Macaulay $\widehat R$-module of dimension $n$. By using [15, Lemma 3.6] and [1, Exercise 8.2.4], we have
$\Gid_{R}\H^{n}_\fm(M)=\Gid_{\widehat R}\H^{n}_{\widehat{\fm}}(M\otimes_R\widehat R)$. Also, in view of [3, Theorem 3.24], $\Gid_{R}M$ and $\Gid_{\widehat R}\widehat M$ are simultaneously finite. Therefore we may assume that $R$
is complete; and hence it has a dualizing complex. Now, one can use 3.8 to obtain the assertion.
\end{proof}
\end{cor}
In [21, Theorem 2.6] a characterization of a complete Gorenstein local ring $R$, in terms of Gorenstein injectivity of the top local cohomology module of $R$, is given. The next corollary together with 2.7 recover that characterization.

\begin{cor}Let $(R,\fm)$ be a Cohen-Macaulay local ring which has a dualizing complex. Then the following conditions are equivalent.
  \begin{itemize}
\item[(i)]{$R$ is Gorenstein.}
\item[(ii)]{$\emph{\Gid}_{R}\emph{\H}^n_{\fa}(R)<\infty$ for any ideal $\fa$ of $R$ such that $R$ is relative Cohen-Macaulay with respect to $\fa$ and that $\emph{\h}_{R}\fa =n$}.
\item[(iii)]{$\emph{\Gid}_{R}\emph{\H}^n_{\fa}(R)<\infty$ for some ideal $\fa$ of $R$ such that $R$ is relative Cohen-Macaulay with respect to $\fa$ and that $\emph{\h}_{R}\fa =n$}.
\end{itemize}
\end{cor}
\begin{proof} The implication (i)$\Rightarrow$(ii) follows from 2.7 and the implication (ii)$\Rightarrow$(iii) is clear. The implication (iii)$\Rightarrow$(i) follows from 3.8(ii) and [3, Proposition 3.11].
\end{proof}

Concerning the above corollary, we notice that if ${\H}^n_{\fa}(R)$ is Artinian, then it is not needed to impose the hypothesis that $R$ has a dualizing complex. Therefore [21, Theorem 2.6] follows from 3.10 without the completeness assumption on $R$.

Next, as promised before, we provide examples to show that if $M$ is not finitely generated or $M$ is not relative Cohen-Macaulay, then 2.5(ii) and 2.5(i), respectively, are no longer true.

\begin{exas}\emph{(i). Let $(R,\fm)$ be a Gorenstein local ring with $\dim R\geq2$ such that $R_{\fp}$ is not regular for some non-maximal prime ideal $\fp$ of $R$ (for example, one can take $R=\frac{K[[X,Y,Z]]}{(X^2)}$ and $\fp=(x,y)R$, where $K$ is a field). Then one can use [3, Theorem 3.14] to see that there exists a non-zero Gorenstein injective $R_{\fp}$-module $M_{\fp}$ which is neither injective nor finitely generated. Hence, by [4, Proposition 10.1.2 ], $\ind_{R_{\fp}} M_{\fp}=\infty$; so that $\ind_{R} M_{\fp}=\infty$. Now, we notice that, for all $x\in\fm-\fp$, $\ind_{R}M_{x}=\infty$ because $(M_{x})_{\fp R_{x}}\cong {M_\fp}$. It is easy to check that $\H^{i}_{\fm}(M_{x})=0$ for all $i$ and that $M_{x}$ is not finitely generated as an $R$-module. Set $N=M_{x}\oplus E_{R}{(k)}$. Then $\ind N=\infty$, but $\Gamma_{\fm}(N)$ is injective. This example shows that, in 2.5(ii), the finiteness assumption on $M$ is required.}

\emph{(ii). Let $R=k[[x,y]]/(xy)$, where $k$ is a field. Then $R$ is a $1$-dimensional complete Gorenstein local ring.  Let $\fm$ be the maximal ideal of $R$ and let $J=(y)R$. In view of [1, Theorem 8.2.1] and [21, Corollary 2.10], $\H_{J}^1(R)$ is a non-zero Gorenstein injective $R$-module. Note that $\Gamma_{J}(R)\neq0$. Now, we show that $\H_{J}^1(R)$ is not injective. If $\H_{J}^1(R)$ were injective, then $\Hom_{R}(\H_{J}^1(R),\E_{R}(R/\fm))= R^n$ for some positive integer $n$.
Therefore, by using [19, Theorem 5.11] and [21, Lemma 3.1], we get an $R$-isomorphism
$$\psi: R^n=\Hom_{R}(\H_{J}^1(R),\E_{R}(R/\fm))\rightarrow \Gamma_{\fm, J}(R)=J.$$ Now, $\psi(x R^n)=xJ=0$ which is a contradiction.
 Hence, by [4, Proposition 10.1.2], we have $\ind{\H_{J}^1(R)}=\infty$. Therefore, by using the exact sequences $$0\longrightarrow\Gamma_{J}(R)\longrightarrow R\longrightarrow R/\Gamma_{J}(R)\longrightarrow0$$ and $$0\longrightarrow R/\Gamma_{J}(R)\longrightarrow R_{y+(xy)}\longrightarrow \H_{J}^1(R)\longrightarrow 0,$$
we achieve $\ind{\Gamma_{J}(R)}=\infty$.}
\end{exas}

As a mentioned just above 3.9, the next theorem is an improvement of [15, Theorem 3.10]. Indeed, we will use the canonical module of a module to prove the above result without assuming that $R$ and $M$ are Cohen-Macaulay. Notice that if $M$ is Cohen-Macaulay, then $K_{M}$ is Cohen-Macaulay. But the converse does not hold in general; see for example [16, Lemma 1.9] and [16, Theorem 1.14].

\begin{thm}Assume that $(R,\fm)$ is local, and $M$ is non-zero finitely generated of dimension $d$. Then the following statements hold.
\begin{itemize}
\item[(i)]{\emph{$\Gid_{R}\H_{\fa}^{d}(M)=\Gpd_{\widehat R}\Gamma_{\fm\widehat R,\fa\widehat R}(K_{\widehat M}).$}
\item[(ii)]{ If $K_{\widehat M}$ is Cohen-Macaulay and \emph{$\Gid_{R}\H_{\fm}^{d}(M)<\infty$}, then
 \emph{$\Gid_{R}\H_{\fm}^{d}(M)=\depth R-d.$}}}
\end{itemize}

\begin{proof}(i) By [1, Theorem 7.1.6], $\H_{\fa}^{d}(M)$ is Artinian. Therefore, by use of [1, Theorem 4.3.2] and [19, Theorem 5.11], we have $\H_{\fa}^{d}(M)=0$ if and only if $\Gamma_{\fm\widehat R,\fa\widehat R}(K_{\widehat M})=0$. Hence, we may assume that $\H_{\fa}^{d}(M)\neq0$. Now, by [15, Lemma 3.6], $\Gid_{R}\H_{\fa}^{d}(M)=\Gid_{\widehat R}\H_{\fa}^{d}(M)$ and, by [3, Theorem 4.25], $\Gid_{\widehat R}\H_{\fa}^{d}(M)=\Gfd_{\widehat R}\Hom_{\widehat R}(\H_{\fa}^{d}(M),\E_{R}(k))$. Therefore, since $\Hom_{\widehat R}(\H_{\fa}^{d}(M),\E_{R}(k))$ is finitely generated as an $\widehat R$-module, one can use [3, Theorem 4.24] and [19, Theorem 5.11] to establish the result.

(ii) First notice that $\Gamma_{\fm\widehat R,\fm\widehat R}(K_{\widehat M})=K_{\widehat M}.$ Hence, by part(i), $\Gid_{R}\H_{\fm}^{d}(M)=\Gpd_{\widehat R}(K_{\widehat M})$. Therefore, in view of [3, Proposition 2.16] and [3, Theorem 1.25], $\Gid_{R}\H_{\fm}^{d}(M)=\depth\widehat R-\depth K_{\widehat M}$. Now, one can use [16, Lemma 1.9(c)] to complete the proof.
\end{proof}
\end{thm}

The following corollary is a generalization of the main result [21, Theorem 2.6].
\begin{cor}Assume that $(R,\fm)$ is local with \emph{$\dim R=d$} and that $K_{\widehat R}$ is Cohen-Macaulay. Then the following statements are equivalent.
\begin{itemize}
\item[(i)]{$R$ is Gorenstein.}
\item[(ii)]{\emph{$\ind_{R}\H^d_{\fm}(R)<\infty$}.}
\item[(iii)]{\emph{$\Gid_{R}\H^d_{\fm}(R)<\infty$}.}
\end{itemize}
\end{cor}
\begin{proof}The implication (i)$\Rightarrow$(ii) follows from 2.5 while the implication (ii)$\Rightarrow$(iii) is clear.

(iii)$\Rightarrow$(i). By 3.12(ii), we have $\Gid_{R}\H_{\fm}^{d}(R)=\depth R-\dim R$; and hence $R$ is Cohen-Macaulay. Now one can use 3.9 to obtain the assertion.
\end{proof}
\begin{cor}Let $(R,\fm)$ be local with $\dim R=d\leq 2$. Then the following statements are equivalent.
\begin{itemize}
\item[(i)]{$R$ is Gorenstein.}
\item[(ii)]{\emph{$\Gid_{R}\H^d_{\fm}(R)<\infty$}.}
\item[(iii)]{\emph{$\H^d_{\fa}(M)$} is Gorenstein injective for all finitely generated $R$--modules $M$ and for all ideals $\fa$ of $R$.}
\end{itemize}
\begin{proof}Let $M$ be a non-zero finitely generated $R$-module. Then, by [1, Theorem 7.1.6], $\H_{\fa}^{d}(M)$ and $\H_{\fm}^{d}(R)$ are Artinian. Therefore, in view of [15, Lemma 3.5], we may assume that $R$ is complete. Since, by [16, Lemma 1.9], $K_{R}$ is Cohen-Macaulay, (i)$\Leftrightarrow$(ii) follows immediately from 3.13. The implication (iii)$\Rightarrow$(i) is clear and the implication (i)$\Rightarrow$(iii) follows from [21, Corollary 2.10].

\end{proof}
\end{cor}
$\mathbf{Acknowledgements}.$ The authors would like to thank Alberto Fernandez Boix and the referee for careful reading of manuscript and helpful comments.



\end{document}